\documentclass[12pt,a4paper]{article}
\usepackage[utf8]{inputenc}
\usepackage[T1]{fontenc}
\usepackage[english]{babel}
\usepackage{amsmath,amssymb,mathtools}
\usepackage{geometry}
\usepackage{graphicx}
\usepackage{subcaption}
\usepackage[
bookmarks=true,
backref=true,
colorlinks=true,
linkcolor=red,
citecolor=blue,
urlcolor=green]{hyperref}
\allowdisplaybreaks

\usepackage{xcolor}
\newcommand{\semi}{\subset \hskip -4mm +}

\numberwithin{equation}{section}

\DeclareMathOperator{\sech}{sech}

\newcommand{\dx}{\partial_x}
\newcommand{\dt}{\partial_t}
\newcommand{\du}{\partial_u}
\newcommand{\dw}{\partial_w}

\newtheorem{remark}{Remark}
\newtheorem{theorem}{Theorem}

\DeclareMathOperator{\cosech}{cosech}

\begin{document}

\title{Lie symmetry analysis of the Saint--Venant equations
with variable bottom}

\author{
Faruk Güngör\textsuperscript{1}\thanks{gungorf@itu.edu.tr},
\quad
Cihangir Özemir\textsuperscript{1}\thanks{Corresponding author, ozemir@itu.edu.tr},
\quad
Pedro J. Torres\textsuperscript{2}\thanks{
ptorres@ugr.es}
\\[1ex]
\small\textsuperscript{1}Department of Mathematics,
Faculty of Science and Letters,\\[-0.2ex]
\small Istanbul Technical University,
34467 Istanbul, Türkiye
\\[0.8ex]
\small\textsuperscript{2}Department of Applied Mathematics,
University of Granada, 18071 Granada, Spain
}

\date{}
\maketitle

\begin{abstract}
The aim of this paper is to apply Lie-symmetry methods to the one-dimensional Saint--Venant shallow-water equations over variable bottom topography. We classify the bottom profiles that admit symmetry extensions, derive similarity reductions to ordinary differential equations, and analyze selected reduced systems. In particular, heteroclinic trajectories yield physically meaningful self-similar solutions describing the draining and filling of water over a quadratic bottom with a fixed shoreline.

\vspace{0.5cm}
\end{abstract}

\section{Introduction}

In this paper, we study the one-dimensional nonlinear shallow-water (Saint--Venant) system over a variable bottom described as
\begin{equation}\label{SV}
\begin{aligned}
&\eta_t+[(\eta+h(x))u]_x=0,\\
&u_t+u u_x+g\eta_x=0,
\end{aligned}
\end{equation}
where $g$ is the gravity constant, $h(x)$ is the bottom profile or undisturbed water depth, $\eta(x,t)$ is the displacement of the free surface from its equilibrium level and $u(x,t)$ is the depth-averaged horizontal fluid velocity. This is the standard one-dimensional model for the evolution of a shallow layer of fluid over a variable topography and has been analyzed by many researchers (see for instance \cite{noble2007,didenkulova2010travelling,Didenkulova_2011,Rybkin_Pelinovsky_Didenkulova_2014,Harris2015} and the references therein).

The gravity constant can be set to $g=1$ by a simple rescaling. Using the total water depth
\begin{equation}\label{SV2}
w(x,t)=\eta(x,t)+h(x),
\end{equation}
as a new dependent variable, system \eqref{SV} is written in the form
\begin{equation}\label{SVw}
\begin{aligned}
&w_t+(uw)_x=0,\\
&u_t+uu_x+w_x=h'(x).
\end{aligned}
\end{equation}
Considering that $w(x,t)$ is the total depth of the water at the point $x$ and the time $t$, $w(x,t)>0$ is the natural admissibility condition for a wet region, while $w(x,t)=0$ describes a shoreline or a dry state.

 The main aim of this paper is to identify novel group-invariant solutions of the system \eqref{SVw} based on the Lie symmetries the system admits for different bottom profiles $h(x)$. Through the identified Lie symmetries, the system can be reduced to a planar system of ordinary differential equations, whose analysis may bring to light some distinguished similarity solutions.

 Lie symmetry analysis is a polyvalent method applied to many problems of Theoretical Physics, and in particular, a Lie symmetry analysis of related shallow water equations can be found in the literature, e.g. \cite{SZATMARI2014530,Pandey2015,Paliathanasis2019,Leach2026}. The symmetry classification problem on uneven bottoms for (2+1) dimensions is studied in \cite{
aksenov2019symmetries,aksenov2020symmetries,bihlo2020lie}. It is necessary to mention that previous works have established the Lie-symmetry classification of the one-dimensional Saint–Venant equations for variable bottom profiles. In \cite{aksenov2016conservation}, the authors determine conservation laws of \eqref{SVw} for certain forms of $h(x)$ and they also introduce the corresponding Lie algebras. They extend their analysis to contact transformations in \cite{aksenov2020conservation}. We mention that quite some time ago a Lie group analysis was successfully  applied to the motion of a rotating shallow liquid in elliptic and circular paraboloidal basins \cite{LeviNucciRogersWinternitz1989} (See also \cite{Chesnokov2009} for a more recent related work).
In the present work, rather than revisiting this classification as an end in itself, we use it to construct optimal systems of subalgebras and systematically derive the associated invariant reductions. We then analyze selected reduced systems qualitatively and explicitly, with particular emphasis on the quadratic-bottom case, where the reduced dynamics possesses heteroclinic connections corresponding to draining and filling shallow-water flows.

The paper is organized as follows. In Section~2, we present the
Lie symmetry algebras of the system \eqref{SVw} for different
bottom profiles $h(x)$. These admissible profiles and the corresponding
Lie symmetry algebras were previously obtained in
Ref.~\cite{aksenov2016conservation}. In the Appendix, we explain how the determining equations yield these allowed forms of the bottom profile. Relying on an existing result of subgroup classification of low dimensional Lie algebras we make use of the optimal system of 1-dimensional subalgebras of these Lie algebras, through which we obtain the reductions of \eqref{SVw} to systems of ODEs in Section 3. In Section 4, we study these systems qualitatively.

\section{The Lie Algebra }

The Lie point symmetries of system (1.3) are generated by
vector fields of the form
\begin{equation}
V=\tau(t,x,w,u)\partial_t
+\xi(t,x,w,u)\partial_x
+\phi(t,x,w,u)\partial_w
+\psi(t,x,w,u)\partial_u.
\end{equation}
We denote the system \eqref{SVw} as $\Delta_1=0$, $\Delta_2=0$.
The infinitesimal invariance criterion requires that the application of the first prolongation $\operatorname{pr}^{(1)}V$ of $V$ to the system be annihilated on the solution manifold:
\begin{equation}
\operatorname{pr}^{(1)}V(\Delta_i)
\big\vert_{\Delta_1=0,\Delta_2=0}=0
\qquad i=1,2.
\end{equation}
After using the original equations to eliminate $w_t$ and $w_x$,
the resulting equations are split with respect to the powers
and products of the remaining independent first-order
derivatives $u_t$ and $u_x$. This yields the determining
equations for $\tau$, $\xi$, $\phi$, and $\psi$.
Eliminating dependent equations, we obtain
the following eight determining equations.

\begin{align}
&h_x\psi+wh_x\psi_w-uh_x\xi_x+u^2h_x\tau_x
 +u\phi_x+w\psi_x-h_x\xi_t+uh_x\tau_t+\phi_t=0,
\label{E1}\\
&\begin{aligned}
\phi-u\psi-uh_x\xi_u+u^2h_x\tau_u
&+u\phi_u+w\psi_u-uwh_x\tau_w-w\phi_w\\
&-uw\psi_w+u^2\xi_x-w\xi_x-u^3\tau_x
+uw\tau_x+u\xi_t-u^2\tau_t+w\tau_t=0,
\end{aligned}
\label{E2}\\
&-\psi-h_x\xi_u+uh_x\tau_u+\phi_u-wh_x\tau_w-w\psi_w
+u\xi_x-u^2\tau_x-w\tau_x+\xi_t-u\tau_t=0,
\label{E4}\\
&\xi_u-u\tau_u+w\tau_w=0,
\label{E6}\\
&-\xi h_{xx}-(h_x)^2\xi_w+u(h_x)^2\tau_w+h_x\phi_w-h_x\xi_x
+uh_x\tau_x+\phi_x+u\psi_x+\psi_t=0,
\label{E7}\\
&\begin{aligned}
\psi-h_x\xi_u+uh_x\tau_u+\phi_u+u\psi_u
&+2uh_x\xi_w-2u^2h_x\tau_w+wh_x\tau_w\\
&-u\phi_w-w\psi_w-u^2\tau_x+w\tau_x-\xi_t=0,
\end{aligned}
\label{E8}\\
&\psi_u+2h_x\xi_w-2uh_x\tau_w-\phi_w
+\xi_x-2u\tau_x-\tau_t=0,
\label{E10}\\
&-\tau_u-\xi_w+u\tau_w=0.
\label{E12}
\end{align}
The determination of the coefficients of the infinitesimal generator \(V\) is a lengthy calculation and is presented in detail in Appendix A. The symmetry generators fall into two main classes, according to whether the bottom profile \(h(x)\) is linear or nonlinear.

The equivalence group of the system is comprised of the transformations
\begin{equation}\label{equiv-tr}
  \tilde{t}=\delta_1 t+\delta_2,  \quad \tilde{x}=\delta_3 x+\delta_4, \quad \tilde{w}=\left(\frac{\delta_3}{\delta_1}\right)^2 w,  \quad \tilde{u}=\frac{\delta_3}{\delta_1} u,
\end{equation}
and $H(x)=h'(x)$ transforms as
$$\tilde{H}(\tilde{x})=\frac{\delta_3}{\delta_1^2}H(\frac{\tilde{x}-\delta_4}{\delta_3}),  \quad \delta_1 \delta_3\ne 0$$ or
$$\tilde{h}(\tilde{x})=\left(\frac{\delta_3}{\delta_1}\right)^2 h(\frac{\tilde{x}-\delta_4}{\delta_3})+C.$$
This group will be specifically useful for detecting discrete symmetries of the system.
\subsection{Symmetry Generators}
\subsubsection{Linear Case, $ h''=0$}
For
\[
h''(x)=0,
\qquad
h(x)=\alpha x+\beta,
\qquad
z=u-\alpha t,
\]
integrating the determining system  will provide us with the following seven infinitesimal generators
\begin{align}
X_1={}&\partial_t,\\
X_2={}&\partial_x,\\
X_3={}&t\partial_x+\partial_u,\\
X_4={}&t\partial_t+\alpha t^2\partial_x
-2w\partial_w+(2\alpha t-u)\partial_u,\\
X_5={}&\left(x-\frac{\alpha}{2}t^2\right)\partial_x
+2w\partial_w+z\partial_u,\\
X_6={}&-3z\,\partial_t
+\left(-\frac{3}{2}u^2+3w+\alpha x+\alpha^2t^2\right)\partial_x
+2\alpha w\,\partial_w-2\alpha z\,\partial_u,\\
X_7={}&\left(x-3tz-\frac{\alpha}{2}t^2\right)\partial_t
+\left(-\frac{3}{2}tu^2+3tw+\alpha tx+\alpha^2t^3\right)\partial_x
+2wz\,\partial_w \notag\\
&+\left(\frac{z^2}{2}+2w+\alpha x
-3\alpha tz-\frac{\alpha^2t^2}{2}\right)\partial_u.
\end{align}
But they do not form a Lie algebra.
This is because the system admits an additional infinite-dimensional family of
hodograph-type symmetry generators
\begin{equation}\label{XFG}
X_{F,G}
=
F(w,z)(\partial_t +\alpha B)+G(w,z)\partial_x,
\end{equation}
where $B=X_3$ is the Galilean boost and $F$ and $G$ satisfy
\begin{equation}\label{eqF}
F_{zz}-wF_{ww}-2F_w=0,
\end{equation}
and
\begin{equation}\label{eqsG}
G_z=zF_z-wF_w,
\qquad
G_w=-F_z+zF_w.
\end{equation}
Eq. \eqref{eqF} is the compatibility condition for the equations in \eqref{eqsG}. The substitution $F=M/w$ removes the last term in \eqref{eqF} so that $M$ satisfies the linear wave equation with variable speed
\begin{equation}
M_{zz}-wM_{ww}=0.
\end{equation}

Now we show that only four of seven generators can span a Lie algebra.
The general infinitesimal generator can be written as
\begin{equation}\label{inflong}
V=t_0X_1+c_3X_2+c_2X_3+t_1X_4
+c_1X_5+p_0X_6+p_1X_7+X_{F,G}.
\end{equation}
This parametrization is not a direct sum in its present form, as
particular choices of \(F\) and \(G\)  reproduce some linear
combinations of \(X_1,\ldots,X_7\).
Indeed, the  generators in our finite list satisfy
\begin{align}
X_1+\alpha X_3&=X_{F=1,\,G=0},\\
X_2&=X_{F=0,\,G=1},\\
X_6-\alpha X_5&=
X_{F=-3z,\,G=-\frac32z^2+3w}.
\end{align}
Thus, $X_1,X_2,X_6$ can be expressed in terms of the
infinite-dimensional family and generators
$X_3,X_5$.  An equivalent parametrization
without this overlap is
\begin{equation}
V=c_3X_3+c_4X_4+c_5X_5+c_7X_7+X_{F,G},
\end{equation}
where $c_3,c_4,c_5,c_7$ are arbitrary constants and
$F,G$ satisfy the auxiliary equations above.

Finally,
let us note that the coefficient of $\partial_\eta$ in the
infinite-dimensional generator printed in Section~3.8 of
\cite{aksenov2016conservation} must be read as
$-\alpha(uw_1+w_2)$ rather than
$+\alpha(uw_1+w_2)$; the negative sign is required by
the invariance criterion.

The nonzero commutation relations for the generators $X_3$, $X_4$, $X_5$ and $X_7$  of the finite dimensional part are
\begin{equation}
\begin{aligned}
[X_3,X_4]&=-X_3, &
[X_3,X_5]&=X_3, &
[X_3,X_7]&=-2X_4-X_5,\\
[X_4,X_7]&=-X_7, &
[X_5,X_7]&=X_7.
\end{aligned}
\end{equation}
The four-dimensional Lie algebra
$\mathfrak{g}_{f}=\langle X_3,X_4,X_5,X_7 \rangle$
is isomorphic to  $ A_{3,8}\oplus A_1 $ in the
Patera--Winternitz classification \cite{patera1977subalgebras}. Indeed, introducing the basis
\begin{equation}
e_1=X_3,\qquad
e_2=-2X_4-X_5,\qquad
e_3=2X_7,\qquad
e_4=X_4+X_5,
\end{equation}
we obtain the nonzero commutation relations
\begin{equation}
[e_1,e_2]=e_1,\qquad
[e_1,e_3]=2e_2,\qquad
[e_2,e_3]=e_3.
\end{equation}
Thus, $\langle e_1,e_2,e_3\rangle$ is isomorphic
to $A_{3,8}\cong\mathfrak{sl}(2,\mathbb{R})$, while $e_4$
commutes with all elements of $\mathfrak{g}_f$ and spans its
one-dimensional center. Consequently,
\begin{equation}
\mathfrak{g}_f
\cong A_{3,8}\oplus A_1
\cong\mathfrak{sl}(2,\mathbb{R})\oplus\mathbb{R}.
\end{equation}
The infinite dimensional part of the Lie algebra is $\mathfrak{g}_{\infty} = \langle X_{F,G}\rangle $.
Introducing the coordinates
\begin{equation}\label{tr-galilei}
y=x-\frac{\alpha}{2}t^2,\qquad z=u-\alpha t,
\end{equation}
and keeping $t$ and $w$ unchanged, the generators take the form
\begin{align}
X_3&=t\partial_y+\partial_z,\\
X_4&=t\partial_t-2w\partial_w-z\partial_z,\\
X_5&=y\partial_y+2w\partial_w+z\partial_z,\\
X_7&=(y-3tz)\partial_t
+t\left(3w-\frac32z^2\right)\partial_y
+2wz\partial_w
+\left(\frac{z^2}{2}+2w\right)\partial_z,\\
X_{F,G}&=F(w,z)\partial_t+G(w,z)\partial_y.
\end{align}
For the commutations of the finite dimensional part with the infinite dimensional part we obtain
\begin{equation}\label{237}
\begin{aligned}
[X_3,X_{F,G}]
&=X_{F_z,\,G_z-F},\\
[X_4,X_{F,G}]
&=X_{-2wF_w-zF_z-F,\,-2wG_w-zG_z},\\
[X_5,X_{F,G}]
&=X_{2wF_w+zF_z,\,2wG_w+zG_z-G},\\
[X_7,X_{F,G}]
&=X_{\widehat F_7,\widehat G_7},
\end{aligned}
\end{equation}
where
\begin{equation}\label{238}
\begin{aligned}
\widehat F_7
&=2wzF_w+\left(\frac{z^2}{2}+2w\right)F_z+3zF-G,\\
\widehat G_7
&=2wzG_w+\left(\frac{z^2}{2}+2w\right)G_z
+\left(\frac32z^2-3w\right)F.
\end{aligned}
\end{equation}
The brackets in \eqref{237}–\eqref{238} preserve the constraints \eqref{eqF}–\eqref{eqsG}.
All functions and their derivatives are evaluated at $(w,z)$. Therefore, we can write
\begin{equation}
\mathfrak{g}(\alpha) =\mathfrak{g}_{\infty} \semi \; \mathfrak{g}_f
               =\langle X_{F,G}\rangle \semi \langle X_3, X_4, X_5, X_7\rangle
\end{equation}
where $\semi$ denotes a semi-direct sum with the abelian ideal on the left.

We comment that the point transformation $(t,x,w,u)\to (t,y,w,z)$ defined by \eqref{tr-galilei}
%\begin{equation}
% y=x-\frac{\alpha}{2}t^2,\quad z=u-\alpha t
%\end{equation}
takes our original system to
\begin{equation}
w_t+(zw)_y=0,\qquad z_t+z z_y+w_y=0.
\end{equation}
Thus,  the constant-slope bottom case is equivalent to the flat-bottom case
in an accelerating coordinate system. The symmetry algebra $\mathfrak{g}(\alpha)$ is isomorphic to $\mathfrak{g}(0)$.

We add that in Ref. \cite{pelinovsky1992exact}, when $h(x)=-x$, the tangential  Legendre transformation was mentioned to transform the system to a Euler-Poisson-Darboux equation (a linear PDE of wave type).

The presence of infinite-dimensional symmetry algebra \eqref{XFG}  suggests that there should exist a hodograph transformation linearising the system. Indeed, this is the case. First of all, we can put $\alpha=0$ (namely, $z=u$) in \eqref{XFG} with no loss of generality. A hodograph transformation is a locally invertible point transformation interchanging the roles of the independent and dependent variables $(t,x,w,u)\to (w,u,t=F,x=G)$ with the non-vanishing Jacobian determinant $J=u_t w_x-u_x w_t\ne 0$. The relations between the first order derivatives of old and new variables are given by the inverse function theorem
\begin{equation}\label{der-rels}
  u_x=J F_w,  \quad u_t=-J G_w,  \quad w_x=-J F_u,  \quad w_t=J G_u.
\end{equation}
Now, it follows that  Equations \eqref{eqsG} are nothing else but the main system \eqref{SVw}, with $h$  linear, written in terms of variables in the hodographic plane $(F,G)$. Again, $F$ satisfies \eqref{eqF}, a linear wave equation with variable speed. We recall that hodograph transformation is closely related to Legendre transformation. This clearly illustrates that the appearance of an infinite-dimensional symmetry in a system of PDEs (or a scalar PDE) can serve as a criterion of  linearisability by  point (or even contact) transformations.

\subsection{Nonlinear Case, $h''\neq 0$}

$V$ is of the form
\begin{equation}\label{R4}
V=(c_0+c_1t)\dt+(c_2x+q(t))\dx+2(c_2-c_1)w\dw
+[(c_2-c_1)u+q'(t)]\du
\end{equation}
with a single restricting equation
\begin{equation}\label{eqh}
[c_2 x+q(t)]h''(x)+(2c_1-c_2)h'(x)-q''(t)=0.
\end{equation}
When we differentiate \eqref{eqh} with respect to $t$ we get $q'(t)h''(x)-q'''(t)=0$, and this gives $q'(t)h'''(x)=0$. Considering the cases $h'''\equiv 0$ and $h'''\not \equiv 0$ separately gives the following admissible profiles of $h(x)$ with the corresponding infinitesimal generators of the associated Lie algebras.

\vspace{12pt}
\noindent\textbf{Arbitrary $h(x)$}. For any function $h(x)$
\begin{equation}
X=\dt
\end{equation}
is a symmetry of \eqref{SVw}.

\begin{itemize}
\item[(I)] $h(x)=h_2x^2+h_1x+h_0$, \quad $h_2>0$.
\begin{align}
X_1&= \partial_t,\\
X_2 &= \big(x+\frac{h_1}{2h_2}\big)\partial_x
      +2w\partial_w+u\partial_u,\\
X_3 &= \exp\big(\sqrt{2h_2}\,t\big)
      \big(\partial_x+\sqrt{2h_2}\,\partial_u\big),\\
X_4 &= \exp\big(-\sqrt{2h_2}\,t\big)
      \big(\partial_x-\sqrt{2h_2}\,\partial_u\big).
\end{align}

We can integrate the symmetry generators as vector fields to find the corresponding symmetry groups. We use them to derive formulas producing new multi-parameter solutions from old solutions.

If the pair $(\widehat U(t,x),\widehat W(t,x))$
constitute a solution of the system, the composite
transformation generates the four-parameter family of solutions
\begin{equation}
\begin{aligned}
 u(t,x)
&=e^{\epsilon_2}\widehat U(\bar t,\bar x)
+\alpha\epsilon_3e^{\alpha t}
-\alpha\epsilon_4e^{-\alpha t},\\
 w(t,x)
&=e^{2\epsilon_2}\widehat W(\bar t,\bar x),
\end{aligned}
\end{equation}
where
\begin{align}
\bar t&=t-\epsilon_1,\\
\bar x&=e^{-\epsilon_2}
\left(x+c-\epsilon_3e^{\alpha t}
-\epsilon_4e^{-\alpha t}\right)-c,
\end{align}
with $\alpha=\sqrt{2h_2}$, $c=h_1/(2h_2)$.

\item[(II)] $h(x)=h_2x^2+h_1x+h_0$, \quad $h_2<0$
\begin{align}
X_1 &= \partial_t,\\
X_2 &= \left(x+\frac{h_1}{2h_2}\right)\partial_x
      +2w\partial_w+u\partial_u,\\
X_3 &= \cos\left(\sqrt{-2h_2}\,t\right)\partial_x
      -\sqrt{-2h_2}\,
      \sin\left(\sqrt{-2h_2}\,t\right)\partial_u,\\
X_4 &= \sin\left(\sqrt{-2h_2}\,t\right)\partial_x
      +\sqrt{-2h_2}\,
      \cos\left(\sqrt{-2h_2}\,t\right)\partial_u.
\end{align}
If $(\widehat U(t,x),\widehat W(t,x))$ are solutions of the system, the composite transformation gives the
four-parameter family of solutions
\begin{align}
 u(t,x)
&=e^{\epsilon_2}\widehat U(\bar t,\bar x)
-\beta\epsilon_3\sin(\beta t)
+\beta\epsilon_4\cos(\beta t),\\
 w(t,x)
&=e^{2\epsilon_2}\widehat W(\bar t,\bar x),
\end{align}
where
\begin{align}
\bar t&=t-\epsilon_1,\\
\bar x&=e^{-\epsilon_2}
\left(x+c-\epsilon_3\cos(\beta t)
-\epsilon_4\sin(\beta t)\right)-c,
\end{align}
with $\beta=\sqrt{-2h_2}$ and $c=h_1/(2h_2)$.

\item[(III)] $h(x)=h_1e^{\lambda x}  + h_0$, \quad $\lambda, h_1\neq 0$.
\begin{align}
X_1 &= \partial_t,\\
X_2 &=   t\dt-\frac{2}{\lambda}\dx-2w\dw-u\du.
\end{align}
Starting from a known solution $(\widehat U(t,x),\widehat W(t,x))$, the Lie group action yields the
two-parameter family of solutions
\begin{align}
 u(t,x)
&=e^{-\epsilon_2}\widehat U\left(
e^{-\epsilon_2}t-\epsilon_1,\,
x+\frac{2\epsilon_2}{\lambda}\right),\\
 w(t,x)
&=e^{-2\epsilon_2}\widehat W\left(
e^{-\epsilon_2}t-\epsilon_1,\,
x+\frac{2\epsilon_2}{\lambda}\right).
\end{align}
\item[(IV)] $h(x)=h_1(x+x_0)^m + h_0$, \quad $m\neq 0,1,2$, $h_1 \neq 0$.
\begin{align}
X_1 &= \partial_t,\\
X_2 &=  (1-\frac{m}{2}\big) t\dt+(x+x_0)\dx+mw\dw+\frac{m}{2}u\du.
\end{align}
The action of the symmetry group on a given solution
$(\widehat U(t,x), \widehat W(t,x))$ produces the
two-parameter family of solutions
\begin{align}
 u(t,x)
&=e^{m\epsilon_2/2}\widehat U\left(
e^{(m/2-1)\epsilon_2}t-\epsilon_1,\,
e^{-\epsilon_2}(x+x_0)-x_0\right),\\
 w(t,x)
&=e^{m\epsilon_2}\widehat W\left(
e^{(m/2-1)\epsilon_2}t-\epsilon_1,\,
e^{-\epsilon_2}(x+x_0)-x_0\right).
\end{align}

\item[(V)] $h(x)=h_1\ln(x+x_0) + h_0$, \quad $h_1\neq 0$.
\begin{align}
X_1 &= \partial_t,\\
X_2 &=   t\dt+(x+x_0)\dx.
\end{align}
For any solution $(\widehat U(t,x),\widehat W(t,x))$,
the symmetry transformations yield the two-parameter
family of solutions
\begin{align}
 u(t,x)
&=\widehat U\left(
e^{-\epsilon_2}t-\epsilon_1,\,
e^{-\epsilon_2}(x+x_0)-x_0\right),\\
 w(t,x)
&=\widehat W\left(
e^{-\epsilon_2}t-\epsilon_1,\,
e^{-\epsilon_2}(x+x_0)-x_0\right).
\end{align}
\end{itemize}
\begin{remark}
The obtained forms of $h(x)$ can further be simplified. In view of Eq. \eqref{SVw}, the constant $h_0$ is of no importance. Moreover, by a translation in the variable $x$, the constant $h_1$ in cases (I) and (II), and the constant $x_0$ in cases (IV) and (V) can be set to zero. However, we prefer to keep those  constants for direct applicability of the results of the manuscript.
\end{remark}
\begin{table}[ht]
\centering
\renewcommand{\arraystretch}{2}
\caption{Commutation relations of the Lie algebras $\mathfrak{g}_I$--$\mathfrak{g}_V$.}
\label{tab:commutation-relations}

\begin{tabular}{|l|ll|}
\hline
\textbf{Algebra}
& \multicolumn{2}{l|}{\textbf{Nonzero Commutation Relations}} \\ \hline

$\mathfrak{g}_I$
& $[X_1,X_3]=\sqrt{2h_2}X_3,$
& $[X_1,X_4]=-\sqrt{2h_2}X_4,$ \\

& $[X_2,X_3]=-X_3,$
& $[X_2,X_4]=-X_4$ \\ \hline

$\mathfrak{g}_{II}$
& $[X_1,X_3]=-\sqrt{-2h_2}X_4,$
& $[X_1,X_4]=\sqrt{-2h_2}X_3,$ \\

& $[X_2,X_3]=-X_3,$
& $[X_2,X_4]=-X_4$ \\ \hline

$\mathfrak{g}_{III}$ & $[X_1,X_2]=X_1$& \\ \hline
$\mathfrak{g}_{IV}$ & $[X_1,X_2]=\big(1-\dfrac{m}{2}\big)X_1$& \\ \hline
$\mathfrak{g}_{V}$ & $[X_1,X_2]=X_1$& \\ \hline
\end{tabular}
\end{table}
Possible forms of $h(x)$ and the  corresponding symmetry generators are available in Reference \cite{aksenov2016conservation}. We include them here as an initial point of our analysis for the study of the reduced equations. As a plus, as \cite{aksenov2016conservation} does not include the solution of the set of overdetermined system of equations leading to possible forms of $h(x)$, we give a detailed treatment of these equations in the Appendix. Now, we present the optimal system of 1-dimensional subalgebras of $\mathfrak{g}_{I},\ldots \mathfrak{g}_{V}$  for a nonlinear $h(x)$.
\subsection{Optimal systems of 1-dimensional subalgebras  for nonlinear  profiles}
\subsubsection{Case I}
By the change of basis
\begin{equation}
e_1=\frac{1}{2}\left(\frac{X_1}{\sqrt{2h_2}}-X_2\right),
\quad e_2=X_3,
\quad
e_3=-\frac{1}{2}\left(\frac{X_1}{\sqrt{2h_2}}+X_2\right),
\qquad e_4=X_4
\end{equation}
the only nonzero commutation relations become
\begin{equation}
[e_1,e_2]=e_2, \qquad [e_3,e_4]=e_4.
\end{equation}
Therefore, in the notation of \cite{patera1977subalgebras}, $\mathfrak{g}_I \simeq 2A_2 = A_2\oplus  A_2 $, that is, the algebra is direct sum of two copies of the non-abelian  two-dimensional Lie algebra. The optimal system of 1-dimensional subalgebras is given in \cite{patera1977subalgebras} as
\begin{equation}\label{opt1}
\langle e_2 \rangle  , \quad \langle e_3 \rangle  , \quad \langle e_4 \rangle  , \quad \langle e_1+ae_3 \rangle  , \quad \langle e_1+\epsilon e_4 \rangle  , \quad \langle e_2+\epsilon e_4 \rangle  , \quad \langle e_2+\epsilon e_3 \rangle
\end{equation}
where $a\in \mathbb{R}$, $\epsilon = \mp 1$.
\subsubsection{Case II}
By the change of basis
\begin{equation}
e_1=X_3,\qquad
e_2=X_4,\qquad
e_3=X_2,\qquad
e_4=-\frac{1}{\sqrt{-2h_2}}X_1
\end{equation}
nonzero commutation relations become
\begin{equation}
[e_1,e_3]=e_1,\qquad
[e_2,e_3]=e_2, \quad
[e_1,e_4]=-e_2,\qquad
[e_2,e_4]=e_1.
\end{equation}
From  \cite{patera1977subalgebras}, we see that  $\mathfrak{g}_{II} \simeq A_{4.12}$. The optimal system of 1-dimensional subalgebras is given in \cite{patera1977subalgebras} as
\begin{equation}\label{opt2}
\langle e_1 \rangle  , \quad \langle e_3 \rangle  , \quad \langle e_4+ae_3 \rangle  .
\end{equation}
with $a\in\mathbb{R}$.

\subsubsection{Cases III, IV, V}
For all of these cases, the Lie algebra is the non-abelian two-dimensional Lie algebra $A_2$ and the optimal system of one-dimensional subalgebras is
\begin{equation}\label{opt345}
\langle X_1 \rangle  , \quad \langle X_2 \rangle  .
\end{equation}
\begin{theorem}
The Lie algebra $\mathfrak{g}$ of the Saint--Venant system \eqref{SVw} is realized depending on the bottom profile $h(x)$ as follows.
\begin{itemize}
\item[(i)] For any bottom profile $h(x)$, the admitted Lie symmetry algebra contains the one-dimensional minimal symmetry algebra
$\mathfrak{g}_0=\langle \dt \rangle  $ corresponding to time translations.
\item[(ii)] For the linear bottom profile, including the flat case, the Lie algebra is infinite dimensional and is introduced in Section 2.1.
\item[(iii)] If the bottom profile is nonlinear, extensions of the minimal symmetry algebra is possible if $h(x)$ is a quadratic polynomial (I,II), is of exponential form (III), is of power form (IV), or
is of logarithmic type (V). In these cases,
\begin{itemize}
\item[(a)] $\dim\mathfrak{g}_{I}=\dim\mathfrak{g}_{II}=4$ and $\mathfrak{g}_{I}$, $\mathfrak{g}_{II}$ are non-isomorphic,
\item[(b)] $\dim\mathfrak{g}_{III}=\dim\mathfrak{g}_{IV}=\dim\mathfrak{g}_{V}=2$ and $\mathfrak{g}_{III}$, $\mathfrak{g}_{IV}$, $\mathfrak{g}_{V}$ are isomorphic.
\end{itemize}
\end{itemize}
\end{theorem}

\section{Reductions to ODEs}
 In this section, we perform reductions of the system \eqref{SVw} to systems of ODEs by making use of the optimal system of 1-dimensional subalgebras obtained above. For each member of the symmetry generators available in \eqref{opt1}, \eqref{opt2} and \eqref{opt345} we will determine the forms of $u(x,t)$ and $w(x,t)$ such that  \eqref{SVw} is reduced to a system of ODEs. The solutions of these ODEs are invariant under the transformations generated by the specific symmetry used in the reduction, therefore they are called group-invariant solutions.

$\dt$ is a symmetry of the system for any $h(x)$. The solutions invariant under this generator are depending on $x$ only. We find them as
\begin{equation}\label{gen1}
u=U(x), \quad w=\frac{W_0}{U(x)}
\end{equation}
where $u(x)$ is found from the solution of the depressed cubic equation
\begin{equation}\label{gen2}
U^3+2(U_0-h(x))U+2W_0=0,
\end{equation}
where $U_0$, $W_0$  are arbitrary constants. This equation is explicitly solvable by Cardano's formula.

\subsection{Case I}
$h(x)=h_2x^2+h_1x+h_0$, \quad $h_2=\alpha^2/2>0$, \quad $\alpha=\sqrt{2h_2}$.

\vspace{12pt}
\noindent\textbf{I.1} Reduction via $\langle e_2 \rangle  $. This case gives a reduction of the form
\begin{equation}
u=\alpha  x +U(t), \quad w=W(t)
\end{equation}
for which we find
\begin{equation}
u=\frac{h_1}{\alpha}+\alpha x +U_0 \exp(-\alpha t), \quad w=W_0 \exp(-\alpha t).
\end{equation}
\noindent\textbf{I.2} Reduction via $\langle e_3 \rangle  $. Solutions are of the form
\begin{equation}
u=\big(x+\frac{h_1}{2h_2}\big)U(\xi), \quad  w=\big(x+\frac{h_1}{2h_2}\big)^2W(\xi), \quad \xi =(x+\frac{h_1}{2h_2})\exp(-\alpha t)
\end{equation}
and we obtain the coupled system
\begin{align}
\xi (\alpha W'-(UW)')-3UW&=0, \\
\xi \big((-\alpha+U)U'+W'\big)+2W+U^2-\alpha^2&=0.
\end{align}
If we set $\rho=\ln|\xi|$, this system becomes
\begin{align}
\alpha W_\rho-(UW)_\rho-3UW&=0, \\
(-\alpha+U)U_\rho+W_\rho+2W+U^2-\alpha^2&=0.
\end{align}

Alternatively, we can search for solutions of the form
\begin{equation}
u=U(\xi)\exp(\alpha t), \quad  w=W(\xi)\exp(2\alpha t), \quad \xi =(x+\frac{h_1}{2h_2})\exp(-\alpha t)
\end{equation}
which results in the equivalent reductions
\begin{align}
\alpha \xi W'-(UW)'-2\alpha W  &=0, \\
 W'-\alpha \xi U'+UU'+\alpha U-\alpha^2\xi&=0.
\end{align}

\noindent\textbf{I.3} Reduction via $\langle e_4 \rangle  $. This case gives a reduction of the form
\begin{equation}
u=-\alpha  x +U(t), \quad w=W(t)
\end{equation}
for which we find
\begin{equation}
u=-\Big(\frac{h_1}{\alpha}+\alpha x\Big) +U_0 \exp(\alpha t), \quad w=W_0 \exp(\alpha t).
\end{equation}
\noindent\textbf{I.4} Reduction via $\langle e_1+ae_3 \rangle  $.

\vspace{12pt}
\noindent\textbf{(i)} If $a=-1$, the symmetry generator is $\dt$. This was treated before in Eqs. \eqref{gen1} and \eqref{gen2}.

\vspace{12pt}
\noindent\textbf{(ii)} When $a=1$, the group-invariant solutions are of the form
\begin{equation}\label{sol1}
u=(x+\frac{h_1}{2h_2})U(t), \quad  w=(x+\frac{h_1}{2h_2})^2W(t)
\end{equation}
We obtain the coupled system of ODEs
\begin{align}
W'+3UW&=0, \label{B2a}\\
U'+U^2+2W&=2h_2. \label{B2b}
\end{align}
We will analyze this system in more detail in Section \ref{sec-4}.

\vspace{12pt}
\noindent\textbf{(iii)} For $a\neq \mp 1$, set $k=\alpha(a+1)/(1-a)$.  Solutions are of the form
\begin{equation}
u=\big(x+\frac{h_1}{2h_2}\big)U(\xi), \quad  w=\big(x+\frac{h_1}{2h_2}\big)^2W(\xi), \quad \xi =\big(x+\frac{h_1}{2h_2}\big)\exp(k t)
\end{equation}
and we obtain the coupled system
\begin{align}
\xi \big[k W'+(UW)'\big]+3UW&=0, \\
\xi \big[(k+U)U'+W'\big]+2W+U^2&=2h_2.
\end{align}
If we set $\rho=\ln|\xi|$, this system becomes
\begin{align}
k W_\rho+(UW)_\rho+3UW&=0, \\
(k+U)U_\rho+W_\rho+2W+U^2&=2h_2.
\end{align}
Using different invariants we can search for alternative solutions of the form
\begin{equation}
u=U(\xi)\exp(-k t), \quad  w=W(\xi)\exp(-2k t), \quad \xi =(x+\frac{h_1}{2h_2})\exp(k t).
\end{equation}
They lead to the reductions
\begin{align}
k \xi W'+(UW)'-2k W  &=0, \\
 W'+k \xi U'+UU'-k  U=2h_2\xi.
\end{align}

\vspace{12pt}
\noindent\textbf{I.5} Reduction via $\langle e_1+\epsilon e_4 \rangle  $.  The solutions are of the form
\begin{equation}
u(x,t)=
\exp(-\alpha t)
\left[U(\xi)-2\epsilon\alpha^2t\right],\quad
w(x,t)
=
\exp(-2\alpha t)W(\xi),
\end{equation}
with
\begin{equation}
\xi
=
\exp(\alpha t)
\left(x+\frac{h_1}{\alpha^2}\right)
-2\epsilon\alpha t
\end{equation}
where $U$, $W$ solve
\begin{align}
(\alpha\xi-2\epsilon\alpha)W'
+(UW)'-2\alpha W&=0,  \label{A51}\\
W'
+\left(U+\alpha\xi-2\epsilon\alpha\right)U'
-\alpha U-2\epsilon\alpha^2-\alpha^2\xi&=0.  \label{A52}
\end{align}
If we define
\begin{equation}
\zeta=\ln (\xi-2\epsilon), \quad U(\xi)=(\xi-2\epsilon) V(\zeta), \quad W(\xi)=(\xi-2\epsilon)^2Z(\zeta)
\end{equation}
the equations \eqref{A51} and \eqref{A52} transform to
\begin{align}
(\alpha+V)\dot Z+Z\dot V+3VZ &=0,\\
\dot Z+\alpha \dot V+V \dot V+2Z+V^2&=\alpha^2+4\epsilon \alpha^2 \exp(-\zeta).
\end{align}
In this form, the system is with constant coefficients and the independent variable only appears in the inhomogeneous term.

%Another possibe reduction is obtained by
%\begin{align}
%u(x,t)&= \left(x+\frac{h_1}{\alpha^2}-2\epsilon \alpha t \exp(-\alpha t)\right) U(\xi)-2\epsilon \alpha^2 t \exp(-\alpha t),\\
%w(x,t)&= \left(x+\frac{h_1}{\alpha^2}-2\epsilon \alpha t \exp(-\alpha t)\right)^2 W(\xi)
%\end{align}
%with
%\begin{equation}
%\xi
%=
%\exp(\alpha t)
%\left(x+\frac{h_1}{\alpha^2}\right)
%-2\epsilon\alpha t
%\end{equation}
%where $U$, $W$ solve
%\begin{align}
%(\alpha\xi^2-2\epsilon\alpha\xi)W'
%+\xi^2(UW)'-4\alpha \epsilon W  +3 \xi UW &=0,\\
%\xi^2W'
%+\left(\xi^2 U+\alpha\xi^2-2\epsilon\alpha\xi\right)U'
%-2\epsilon \alpha U+2\xi W+\xi U^2-2\epsilon\alpha^2-\alpha^2\xi&=0.
%\end{align}

\vspace{12pt}
\noindent\textbf{I.6} Reduction via $\langle e_2+\epsilon e_4 \rangle  $.

\vspace{12pt}
For $\epsilon=1$, the group-invariant solutions are of the form
\begin{align}
u(x,t)
&=
\alpha\left(x+\frac{h_1}{\alpha^2}\right)
\tanh(\alpha t)+U(t),\\
w(x,t)&=W(t).
\end{align}
The reduced system is
\begin{align}
W'+\alpha\tanh(\alpha t)W&=0,\\
U'+\alpha\tanh(\alpha t)U&=0.
\end{align}
Hence,
\begin{align}
u(x,t)
&=
\alpha\left(x+\frac{h_1}{\alpha^2}\right)
\tanh(\alpha t)\label{sol-u}
+U_0\operatorname{sech}(\alpha t),\\
w(x,t)&=W_0\operatorname{sech}(\alpha t).\label{sol-w}
\end{align}
For $\epsilon=-1$, the group-invariant solutions are of the form
\begin{align}
u(x,t)
&=
\alpha\left(x+\frac{h_1}{\alpha^2}\right)
\coth(\alpha t)+U(t),\\
w(x,t)&=W(t).
\end{align}
The reduced system is
\begin{align}
W'+\alpha\coth(\alpha t)W&=0,\\
U'+\alpha\coth(\alpha t)U&=0.
\end{align}
Therefore,
\begin{align}
u(x,t)
&=
\alpha\left(x+\frac{h_1}{\alpha^2}\right)
\coth(\alpha t)
+U_0\operatorname{csch}(\alpha t),\label{sol-u2}\\
w(x,t)&=W_0\operatorname{csch}(\alpha t).\label{sol-w2}
\end{align}

\vspace{12pt}
\noindent\textbf{I.7} Reduction via $\langle e_2+\epsilon e_3 \rangle  $. The group-invariant solutions are of the form
\begin{align}
u(x,t)
&=\exp(\alpha t)
\left[U(\xi)-2\epsilon\alpha^2t\right],\\
w(x,t)
&=\exp(2\alpha t)W(\xi),\\
\xi
&=\exp(-\alpha t)
\left(x+\frac{h_1}{\alpha^2}\right)
+2\epsilon\alpha t.
\end{align}
The corresponding reduced system is
\begin{align}
(2\epsilon\alpha-\alpha\xi)W'
+(UW)'+2\alpha W&=0,\\
W'
+\left(U-\alpha\xi+2\epsilon\alpha\right)U'
+\alpha U-2\epsilon\alpha^2-\alpha^2\xi&=0.
\end{align}
\begin{remark}
Although the reductions associated with
$\langle e_1+\epsilon e_4\rangle$ and
$\langle e_2+\epsilon e_3\rangle$ appear as distinct cases in the
optimal system constructed under the adjoint action of the connected
Lie group, they are related by the discrete symmetry
\begin{equation}
(t,x,u,w)\mapsto(-t,x,-u,w).
\end{equation}
of the system \eqref{SVw}.
Equivalently, their formulas are interchanged formally by
$\alpha\mapsto-\alpha$; however, since
$\alpha=\sqrt{2h_2}>0$ has been fixed, the discrete symmetry provides
the appropriate equivalence for fixed $\alpha$.\end{remark}
\subsection{Case II}
$h(x)=h_2x^2+h_1x+h_0$, \quad $h_2=-\beta^2/2<0$, \quad $\beta=\sqrt{-2h_2}$.

\vspace{12pt}
\noindent\textbf{II.1} Reduction via $\langle e_1 \rangle  $. This case gives a reduction of the form
\begin{equation}
u=-\beta x \tan(\beta t)+U(t), \quad w=W(t)
\end{equation}
for which the resulting system of ODE is solvable. We find
\begin{equation}\label{sing}
u=(\frac{h_1}{\beta}-\beta x) \tan(\beta t)+U_0\sec{\beta t}, \quad w=W_0 \sec(\beta t).
\end{equation}
\noindent\textbf{II.2} Reduction via $\langle e_3 \rangle  $. In this case, the group-invariant solutions are of the form
\begin{equation}
u=(x+\frac{h_1}{2h_2})U(t), \quad  w=(x+\frac{h_1}{2h_2})^2W(t)
\end{equation}
which was studied above in \textbf{I.4 (ii)}.

\vspace{12pt}
\noindent\textbf{II.3} Reduction through $\langle e_4+ae_3 \rangle  $. Solutions are of the form
\begin{equation}
u=(x+\frac{h_1}{2h_2})U(\xi), \quad  w=(x+\frac{h_1}{2h_2})^2W(\xi), \quad \xi =(x+\frac{h_1}{2h_2})\exp(a\beta t)
\end{equation}
and we obtain the coupled system
\begin{align}
\xi (a\beta W'+(UW)')+3UW&=0, \\
\xi \big((a\beta+U)U'+W'\big)+2W+U^2+\beta^2&=0.
\end{align}
For a different selection of the invariants, with
\begin{equation}
u=U(\xi)\exp(-a\beta t), \quad  w=W(\xi)\exp(-2a\beta t), \quad \xi =(x+\frac{h_1}{2h_2})\exp(a\beta t)
\end{equation}
we find
\begin{align}
a\beta \xi W'+(UW)'-2a\beta W  &=0, \\
 W'+a\beta \xi U'+UU'-a\beta  U+\beta^2\xi&=0.
\end{align}

\subsection{Cases III, IV, V}
In each case, the optimal system of one-dimensional subalgebras includes  $X_1=\langle \dt \rangle  $, for which the corresponding reduced equation was derived above. The remaining representative $X_2$ generates a scaling transformation in each case. We present these reductions separately below.

\vspace{12pt}
\noindent\textbf{Case III} $h(x)=h_1 e^{\lambda x}+h_0$

\vspace{12pt}
\noindent   Group-invariant solutions have the form
\begin{equation}
u(x,t)=e^{\lambda x/2} U(\xi), \quad w(x,t)=e^{\lambda x} W(\xi), \quad \xi=te^{\lambda x/2}
\end{equation}
yielding the reduced system
\begin{align}
2W'+ \lambda\xi  (UW)'+3\lambda UW  &=0,\\
2U'+ \lambda\xi  U U'+\lambda\xi  W'+2 \lambda W+\lambda U^2&=2 h_1 \lambda.
\end{align}

\noindent\textbf{Case IV} $h(x)=h_1 x^m+h_0$, $m\neq 0,1,2$

\vspace{12pt}
\noindent We have removed $x_0$ by a translation in $x$. With
\begin{equation}
u(x,t)=x^{m/2} U(\xi), \quad w(x,t)=x^{m} W(\xi), \quad \xi=tx^{\frac{m}{2}-1}
\end{equation}
we obtain
\begin{align}
2W'+(m-2)\xi UW'
+(m-2)\xi W U'
+3mUW&=0,\\
2U'+(m-2)\xi UU'
+(m-2)\xi W'
+mU^2+2mW&=2mh_1.
\end{align}

\noindent\textbf{Case V} $h(x)=h_1 \ln (x+x_0)+h_0$

\vspace{12pt}
\noindent For this last case we have
\begin{equation}
u(x,t)=U(\xi), \quad w(x,t)=W(\xi), \quad \xi=\dfrac{x+x_0}{t}
\end{equation}
with the reduced system
\begin{align}
\xi  W'-(UW)'&=0,\\
\xi U'-UU'-W'+\frac{h_1}{\xi}&=0.
\end{align}

\section{Analysis of the reduced equations and physical interpretation of the solutions}\label{sec-4}

\subsection{Spatially uniform total-depth solutions}
In the analysis presented in Section 3, we have obtained the explicit expression of some solutions for which the total depth $w$ is independent of $x$, although it may vary with time. We call this type of solutions {\it spatially uniform total-depth solutions}.

\begin{figure}[ht]
    \centering
\begin{subfigure}{0.48\textwidth}
        \centering
        \includegraphics[width=\linewidth]{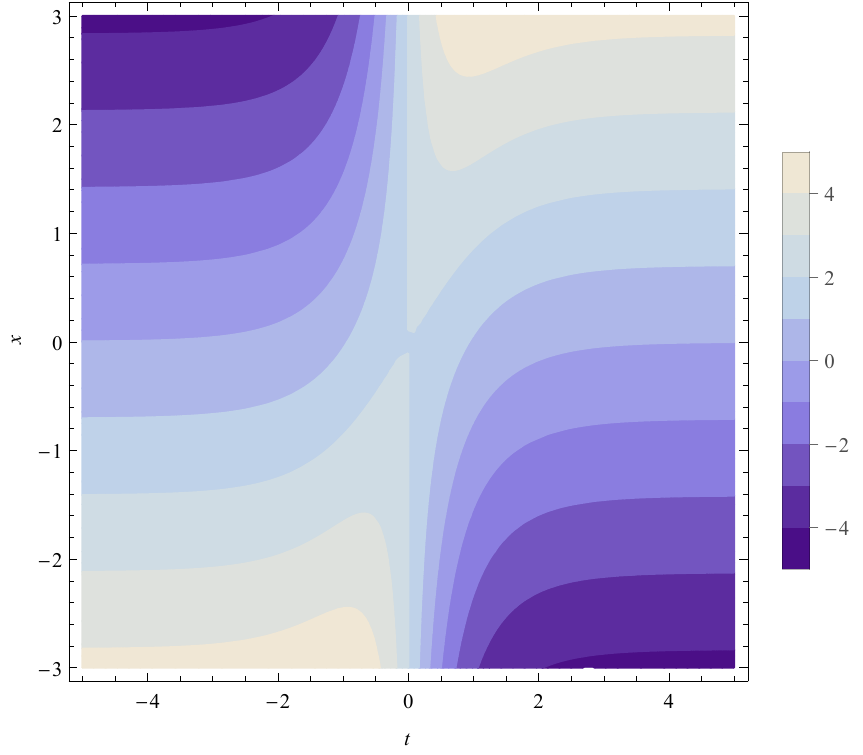}
        \caption{Horizontal velocity $u(x,t)$}
        \label{fig:figura1}
    \end{subfigure}
    \hfill
    \begin{subfigure}{0.48\textwidth}
        \centering
        \includegraphics[width=\linewidth]{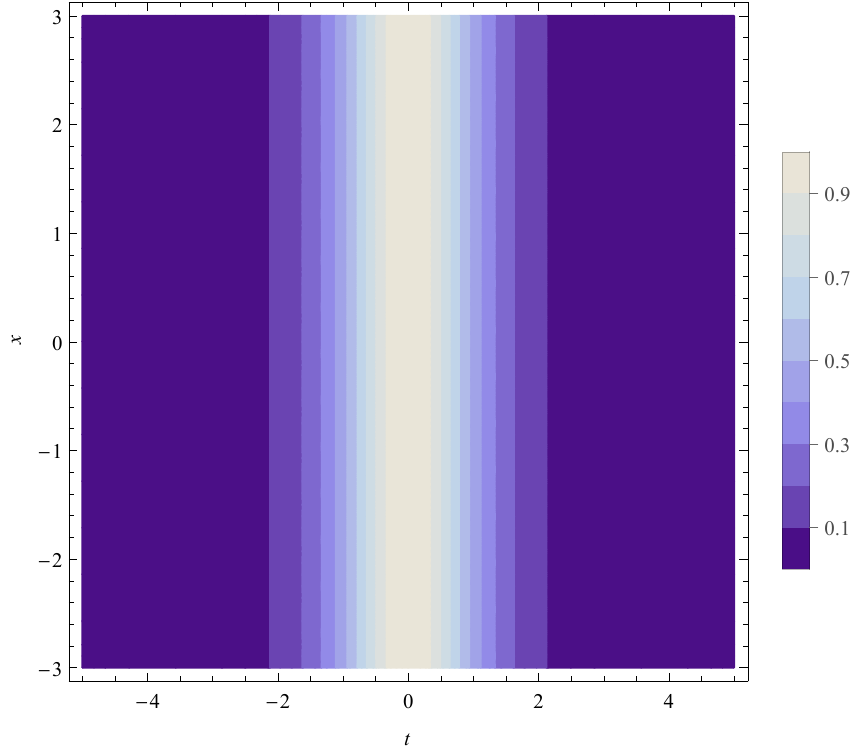}
        \caption{Total water depth $w(x,t)$}
        \label{fig:figura2}
    \end{subfigure}

    \caption{The time-localized depth pulse \eqref{sol-u}-\eqref{sol-w} for $h_1=0$, $h_2=1$, $U_0=2,W_0=1$.}
    \label{fig1}
\end{figure}

In particular, for the case of a quadratic bottom (Case I) $h(x)=h_2x^2+h_1x+h_0$ with $h_2>0$, we obtained in \eqref{sol-u}-\eqref{sol-w} the explicit solution
\begin{align*}
u(x,t)
&=
\alpha\left(x+\frac{h_1}{\alpha^2}\right)
\tanh(\alpha t)
+U_0\operatorname{sech}(\alpha t),\\
w(x,t)&=W_0\operatorname{sech}(\alpha t),
\end{align*}
where $\alpha=\sqrt{2h_2}$, and $U_0,W_0$ are arbitrary constants. As $w$ is the total water depth, only the case $W_0>0$ is physically significant. This solution is identified as a two parametric family of time-localized depth pulses, i.e., spatially uniform total-depths solutions that are localized in time. The positive water depth depends only on $t$ and tends to zero as $t\to\pm\infty$. Therefore, these solutions connect two asymptotically dry states. Fig. \ref{fig1} presents the contour plots of an example of solution of this family.

Other explicit solutions obtained, such as \eqref{sol-u2}-\eqref{sol-w2} or \eqref{sing}, present singularities and have a less clear physical interpretation.

It is interesting to note that the spatially uniform total-depth solutions can be derived without using a Lie symmetry analysis. In order to see this, let us assume $w=W(t)$. Then  the first equation of the system  \eqref{SVw} is integrated to
\begin{equation}\label{u-lin}
  u(x,t)=-\frac{\dot{W}(t)}{W(t)}x+U(t).
\end{equation}
Substituting $u$ in the second equation of the system with quadratic bottom $h(x)=h_2x^2+h_1 x+h_0$  we find a linear polynomial in $x$ with time dependent coefficients. Setting the coefficients of $x$ and $1$ to zero we get a second order nonlinear ODE for $W(t)$ and a linear ODE for $U(t)$
\begin{equation}\label{ode-W-U}
 W\ddot{W}-2\dot{W}^2+\alpha^2 W^2=0, \quad \dot{U}-\frac{\dot{W}(t)}{W(t)} U=h_1,  \quad \alpha^2=2h_2.
\end{equation}
The equation for $W$ is equivalent to  the linear ODE $\ddot{Z}-\alpha^2 Z=0$ by the reciprocal change of the dependent variable $W=Z^{-1}$. So we find $W(t)=W_0 \sech(\alpha t+W_1)$ or $W(t)=W_0 \cosech(\alpha t+W_1)$, where $W_0, W_1$ are integration constants. Then $U(t)$ is easy to obtain
$$U(t)=\left(h_1\int W^{-1}(t)dt+U_0\right)W=\frac{h_1}{\alpha}\tanh (\alpha t+W_1)+U_0 W_0\sech (\alpha t+W_1),$$
where $U_0$ is a constant. The solution of the system is then
\begin{equation}\label{sol-gen}
  u=\left(\alpha x+\frac{h_1}{\alpha}\right)\tanh (\alpha t+W_1)+U_0W_0\sech (\alpha t+W_1),  \quad w=W_0 \sech (\alpha t+W_1).
\end{equation}

If $2h_2=-\beta^2$, then we replace $\alpha^2\to -\beta^2$ and the hyperbolic functions by trigonometric ones. As a result, when $w$ merely depends on time, all invariant solutions obtained in Cases I and II can be recovered as specific cases of these general solutions.

\subsection{Qualitative phase plane analysis of system \eqref{B2a}-\eqref{B2b}}\label{subsec4.2}

Among all the reduced systems obtained, only the analysis of system \eqref{B2a}-\eqref{B2b} is easy to handle. We provide two complementary approaches. In this subsection, we will perform a phase-plane analysis, while in the next subsection some explicit solutions are obtained.

\begin{figure}[htbp]
    \centering
\includegraphics[width=12cm]{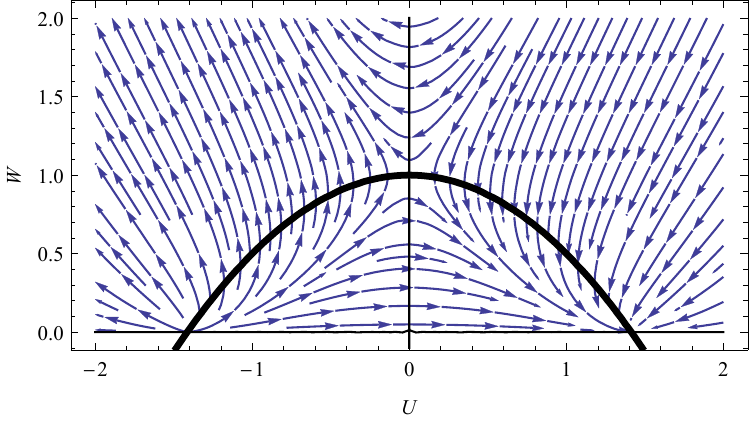}
\caption{Vector field of system \eqref{B2a}-\eqref{B2b} with $h_2=1$, where the depicted parabola is the vertical nullcline.}
\label{fig2}       % Give a unique label
\end{figure}

 For any $h_2>0$, there are exactly three equilibrium points
\[
E_1=(0,h_2),\quad E_2=(\sqrt{2h_2},0),\quad E_3=(-\sqrt{2h_2},0).
\]
An elementary linear analysis shows that $E_1$ is a saddle, $E_2$ is a stable node, and $E_3$ is an unstable node. The horizontal nullclines are the axes $U=0,W=0$, while the vertical nullcline is the parabola $W=h_2-U^2/2$. Nullclines separate the phase plane into regions where the components of the vector field have constant sign, as shown in Fig.\ref{fig2}. Note that from \eqref{B2a} the horizontal axis $W=0$ is invariant, therefore any solution $(U,W)$ has a constant sign in $W$.

From an elementary analysis of the phase plane, one deduces that there exists
\begin{itemize}
    \item one heteroclinic orbit connecting $E_1$ to $E_2$
     \item one heteroclinic orbit connecting $E_3$ to $E_1$
      \item infinitely many heteroclinic orbits connecting $E_3$ to $E_2$.
\end{itemize}

All these heteroclinic connections belong to the semiplane $W>0$ and hence have a physical meaning.The heteroclinic orbit connecting $E_1$ to $E_2$ corresponds to the branch of the unstable manifold of $E_1$ in the region delimited by the vertical nullcline (the parabola depicted in Fig.\ref{fig2}) and the horizontal axis. It would correspond to a solution $(U,W)$ of system \eqref{B2a}-\eqref{B2b} such that
\begin{equation}
\lim_{t\to-\infty}U(t)=0,\lim_{t\to+\infty}U(t)=\sqrt{2h_2},\lim_{t\to-\infty}W(t)=h_2,\lim_{t\to+\infty}W(t)=0.
\end{equation}
Moreover, both components $U$ and $W$ are positive for any $t$. In view of \eqref{sol1}, the original Saint--Venant system \eqref{SVw} has a solution
$$
u=(x+\frac{h_1}{2h_2})U(t), \quad  w=(x+\frac{h_1}{2h_2})^2W(t),
$$
which is clearly interpretable physically: the total water depth $w(x,t)$ is zero if and only if $x=-\frac{h_1}{2h_2}$, which is the vertex of the quadratic profile $h(x)$ and can be interpreted as a stationary dry point or fixed shoreline, and tends to zero when $t\to+\infty$ for any fixed $x$. We call a solution with this property a {\it draining solution}. The interpretation of the wave behaviour when $t\to -\infty$ is also clear: the total water depth $w(x,t)$ tends to the equilibrium profile $(x+\frac{h_1}{2h_2})^2h_2$.

The infinite number of heteroclinic orbits connecting $E_3$ to $E_2$ have an analogous physical interpretation: they can be regarded as time-localized waves or positive-depth pulses connecting two asymptotically dry states, in analogy with the explicit solution  \eqref{sol-u}-\eqref{sol-w}. However, these solutions are different from \eqref{sol-u}-\eqref{sol-w} because of the explicit dependence of the water depth $w$ on $x$ and the presence of the fixed shoreline at $x=-\frac{h_1}{2h_2}$.

Similarly, the heteroclinic orbit connecting $E_3$ to $E_1$ corresponds to the branch of the stable manifold of $E_1$ in the region delimited by the vertical nullcline and the horizontal axis.  In this case, the solution $(U,W)$ of system \eqref{B2a}-\eqref{B2b} verifies
\begin{equation}
\lim_{t\to-\infty}U(t)=-\sqrt{2h_2},\lim_{t\to+\infty}U(t)=0,\lim_{t\to-\infty}W(t)=0,\lim_{t\to+\infty}W(t)=h_2.
\end{equation}
 The component $U$ is negative while $W$ is positive for any $t$, and both are increasing functions of $t$. The corresponding solution of the Saint--Venant system increases the total wave depth from $w=0$ at $t\to-\infty$ to the fixed equilibrium profile $w=(x+\frac{h_1}{2h_2})^2h_2$ when $t\to+\infty$. We call it a {\it filling solution}.

\subsection{Explicit solutions  of system \eqref{B2a}-\eqref{B2b}}

If we solve $W$ from \eqref{B2b} as
\begin{equation}
W=(2h_2-U'-U^2)/2
\end{equation}
and substitute in \eqref{B2a}, we obtain
\begin{equation}\label{reduced}
U''+5UU'-6h_2U+3U^3=0.
\end{equation}

We consider the second-order ordinary differential equation \eqref{reduced}.
 Since the independent variable $t$
does not appear explicitly ($\dt$ invariance), the order of the equation can be reduced
by using the differential invariants
$U$, $p(U)=U'(t)$ and $p'(U)=U''/U'=U''/p$.
In terms of the invariants, \eqref{reduced} takes the form
\begin{equation}\label{eq:Abel_U}
p\frac{dp}{dU}+5Up+3U(U^2-2h_2)=0,
\end{equation}
which is an Abel equation of the second kind. The constant solutions,
which must be considered separately in this reduction, satisfy
\[
3U(U^2-2h_2)=0,
\]
and corresponds to the equilibria of system \eqref{B2a}-\eqref{B2b} yet identified in Subsection \ref{subsec4.2}.

For the particular coefficients in \eqref{eq:Abel_U}, a further
substitution transforms the equation into a separable form. Setting
\begin{equation}\label{eq:Abel_substitution}
p=(U^2-2h_2)v(U)
\end{equation}
and working on an interval where $U^2-2h_2\neq0$, we obtain
\[
\frac{dp}{dU}=2Uv+(U^2-2h_2)\frac{dv}{dU}.
\]
Substitution into \eqref{eq:Abel_U} yields
\begin{equation}\label{eq:Abel_separable}
(U^2-2h_2)v\frac{dv}{dU}
+U(v+1)(2v+3)=0.
\end{equation}
After separating variables, integrating and then returning to $p$ through \eqref{eq:Abel_substitution}, we obtain
\[
\frac{\bigl[2p+3(U^2-2h_2)\bigr]^3}
     {\bigl[p+U^2-2h_2\bigr]^2}=C_1.
\]
Thus, a first integral of \eqref{reduced} is
\begin{equation}\label{eq:U_first_integral}
\frac{\bigl[2U'+3(U^2-2h_2)\bigr]^3}
     {\bigl[U'+U^2-2h_2\bigr]^2}=C_1,
\end{equation}
where the denominator is assumed to be nonzero. Equivalently,
\begin{equation}\label{eq:U_algebraic_integral}
\bigl[2U'+3(U^2-2h_2)\bigr]^3
=C_1\bigl[U'+U^2-2h_2\bigr]^2.
\end{equation}
For a fixed value of $C_1$, this relation defines $U'=p(U)$
implicitly through a cubic algebraic equation. Once a real branch
$p(U)$ has been selected, the corresponding solution is determined
locally by the quadrature
\begin{equation}\label{eq:U_quadrature}
t-t_0=\int_{U_*}^{U(t)}\frac{ds}{p(s)},
\end{equation}
on an interval where $p(s)\neq0$.

The values $v=-1$ and $v=-3/2$, excluded when separating
\eqref{eq:Abel_separable}, give two explicit solution families.
Indeed, they lead respectively to
\begin{equation}\label{eq:U_special_branches}
U'=2h_2-U^2,
\qquad
U'=\frac32(2h_2-U^2).
\end{equation}
The first of these branches makes the denominator in
\eqref{eq:U_first_integral} vanish and must therefore be retained
separately. The second corresponds to $C_1=0$ away from
$U^2=2h_2$.

Since $h_2>0$ and $\alpha =\sqrt{2h_2}$, both equations in \eqref{eq:U_special_branches} have the form
\[
U'=\kappa(\alpha^2-U^2),
\qquad
\kappa\in\left\{1,\frac32\right\}.
\]
For nonconstant solutions, integration gives
\[
\frac{1}{2\alpha}\ln\left|\frac{\alpha+U}{\alpha-U}\right|
=\kappa(t-t_0).
\]
For $|U|<\alpha$, this yields the smooth solutions
\begin{equation}\label{eq:U_tanh_solutions}
U(t)=\alpha\tanh\bigl(\alpha(t-t_0)\bigr),
\qquad
U(t)=\alpha\tanh\left(\frac{3\alpha}{2}(t-t_0)\right).
\end{equation}
These solutions are defined for all $t\in\mathbb{R}$ and satisfy
\[
\lim_{t\to-\infty}U(t)=-\alpha,
\qquad
\lim_{t\to+\infty}U(t)=\alpha.
\]
For $\kappa=1$, the couple of solutions $u$, $w$ is
\begin{equation}
u(x,t)=\alpha \big(x+\frac{h_1}{\alpha^2}\big)\tanh\bigl(\alpha(t-t_0)\bigr), \quad w(x,t)=0,
\end{equation}
and for  $\kappa=3/2$ we have
\begin{align}
u(x,t)&=\alpha \big(x+\frac{h_1}{\alpha^2}\big)\tanh\left(\frac{3\alpha}{2}(t-t_0)\right),\\
w(x,t)&=-\frac{\alpha^2}{4}\big(x+\frac{h_1}{\alpha^2}\big)^2 \mathrm{sech}^2 \left(\frac{3\alpha}{2}(t-t_0)\right).
\end{align}
This solution is mathematically valid for the reduced ODEs, but it is not physically admissible as a water-depth solution because $w<0$.

\begin{figure}[ht]
    \centering
\begin{subfigure}{0.48\textwidth}
        \centering
        \includegraphics[width=\linewidth]{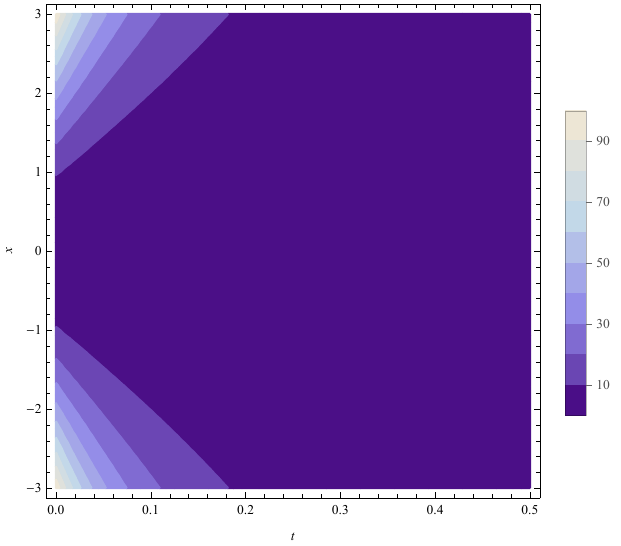}
        \caption{Horizontal velocity $u(x,t)$}

    \end{subfigure}
    \hfill
    \begin{subfigure}{0.48\textwidth}
        \centering
        \includegraphics[width=\linewidth]{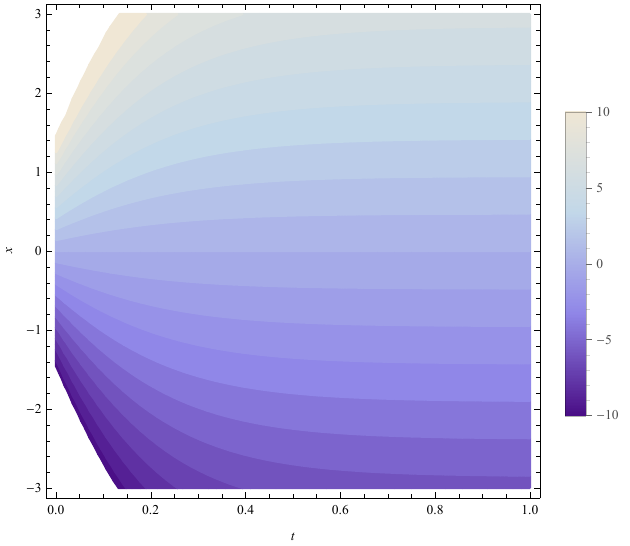}
        \caption{Total water depth $w(x,t)$}
        \label{fig:figura2}
    \end{subfigure}
     \caption{The solution \eqref{4.19}-\eqref{4.20} for $h_1=0$, $h_2=2$ and $t_0=-0.1$.}
    \label{fig3}
\end{figure}

For $|U|>\alpha$, integration instead gives the singular solutions
\begin{equation}\label{eq:U_coth_solutions}
U(t)=\alpha\coth\bigl(\alpha(t-t_0)\bigr),
\qquad
U(t)=\alpha\coth\left(\frac{3\alpha}{2}(t-t_0)\right).
\end{equation}
Each of these solutions has a pole at $t=t_0$ and is defined
separately on the intervals $(-\infty,t_0)$ and $(t_0,\infty)$.
In this case, for $\kappa=1$, the couple of solutions $u,w$ is
\begin{equation}
u(x,t)=\alpha \big(x+\frac{h_1}{\alpha^2}\big)\coth\bigl(\alpha(t-t_0)\bigr), \quad w(x,t)=0,
\end{equation}
which is a dry state with no physical interest, while for  $\kappa=3/2$ we have
\begin{align}
u(x,t)&=\alpha \big(x+\frac{h_1}{\alpha^2}\big)\coth\left(\frac{3\alpha}{2}(t-t_0)\right),\label{4.19}\\
w(x,t)&=\frac{\alpha^2}{4}\big(x+\frac{h_1}{\alpha^2}\big)^2 \mathrm{cosech}^2 \left(\frac{3\alpha}{2}(t-t_0)\right).\label{4.20}
\end{align}
We depict in Fig. \ref{fig3} the contour plots of one of the solutions of this family.

\section{Conclusion}

In this paper, we have investigated the one-dimensional Saint–Venant equations over variable bottom topography by means of Lie symmetry analysis. Using the total water depth $w$ and the depth-averaged horizontal velocity $u$ as dependent variables, we have examined the symmetry algebras admitted by different bottom profiles, constructed optimal systems of one-dimensional subalgebras for the nonlinear cases, and derived the corresponding reductions to ordinary differential equations. Although the symmetry classification itself has been established in a previous work \cite{aksenov2016conservation}, our analysis provides a systematic framework and includes a detailed account of solving the determining system of equations for the infinitesimal coefficients in the Appendix.

The symmetry structure depends strongly on the bottom topography. For a linear bottom profile, the symmetry algebra is the semi-direct sum of a four-dimensional Lie algebra with an infinite-dimensional Lie algebra. Further, in this case, the system is equivalent to the flat-bottom flow. Among the nonlinear bottom profiles, only quadratic, power-law, exponential and logarithmic topographies allow nontrivial Lie symmetries. The two quadratic cases, one up and one down, admit the richest symmetry structure: Maximal dimension of Lie symmetry algebra of the Saint--Venant system  with a nonlinear bottom profile is four and is realised for parabolic bottom profiles.  For nonlinear bottom profiles, we have presented the optimal system of one-dimensional subalgebras of the related Lie algebra. Through these symmetry generators, we obtained reductions of the original shallow-water equations to planar systems of ODEs.

These results show how the geometry of the bottom constrains the transformations that preserve the shallow-water equations and, consequently, the types of similarity reductions available. However, symmetry classification should be distinguished from the physical classification of solutions: the existence of a symmetry or an exact solution does not guarantee that the resulting water depth is positive or that the solution is free from singularities by itself. For this reason, we have paid special attention to the qualitative behavior and explicit integration of the reduced ODE systems, leading to non-trivial similarity solutions with distinct physical mechanisms. In this way, we have completely described the spatially uniform depth solutions. In contrast, the solutions coming from the qualitative analysis of the reduced system \eqref{B2a}-\eqref{B2b} have an explicit spatial dependence in the total depth and retain a fixed shoreline associated with the vertex of the quadratic bottom. We have identified three equilibria and several types of heteroclinic connections on the reduced dynamics. These trajectories describe transitions between different asymptotic states and, through the similarity ansatz, generate novel nontrivial solutions identified as draining or filling flows of the original shallow water equations with quadratic bottom.

\section*{Acknowledgements}
C. Özemir visited the University of Granada during the preparation of this manuscript and would like to express his sincere gratitude to Prof. P. J. Torres and the members of the Department for their warm hospitality.

\section*{Funding}
This work was supported by Scientific Research Projects Department of Istanbul Technical University, Project Number: TGA-2026-48385.

%\bibliographystyle{unsrt}
%\bibliography{NLVCWE}

\newpage

\appendix
\section{Appendix: Solution of the Determining Equations}

From \eqref{E6} and \eqref{E12},
\begin{align}
\xi_u&=u\tau_u-w\tau_w,\label{xi-u}\\
\xi_w&=-\tau_u+u\tau_w\label{xi-w}
\end{align}
are obtained. The derivatives of \eqref{xi-u} with respect to $w$ and of
\eqref{xi-w} with respect to $u$ are
\begin{align}
\xi_{uw}&=u\tau_{uw}-\tau_w-w\tau_{ww},\\
\xi_{wu}&=-\tau_{uu}+\tau_w+u\tau_{uw}.
\end{align}
From the compatibility condition $\xi_{uw}=\xi_{wu}$
\begin{equation}
\tau_{uu}-w\tau_{ww}-2\tau_w=0
\label{T1}
\end{equation}
is obtained.

If \eqref{E6} is used in \eqref{E4}, the first group containing $h_x$ disappears and
\begin{equation}
\phi_u-w\psi_w-\psi+u\xi_x-(u^2+w)\tau_x+\xi_t-u\tau_t=0
\label{E4s}
\end{equation}
is found. Similarly, if \eqref{E12} is used in \eqref{E10},
\begin{equation}
\psi_u-\phi_w+\xi_x-2u\tau_x-\tau_t-2h_x\tau_u=0
\label{E10s}
\end{equation}
is obtained. From here
\begin{align}
\phi_u&=\psi+w\psi_w-u\xi_x+(u^2+w)\tau_x-\xi_t+u\tau_t,
\label{phi-u}\\
\phi_w&=\psi_u+\xi_x-2u\tau_x-\tau_t-2h_x\tau_u.
\label{phi-w}
\end{align}
Now \eqref{phi-u} and \eqref{phi-w} are used in \eqref{E8}. Additionally
$\xi_u$, $\xi_w$ terms are eliminated with the help of \eqref{xi-u} and \eqref{xi-w}. After simplification, one finds
\begin{equation}
\psi-u\xi_x+(u^2+w)\tau_x-\xi_t+u\tau_t+wh_x\tau_w=0
\end{equation}
and therefore
\begin{equation}
\psi=u\xi_x-(u^2+w)\tau_x+\xi_t-u\tau_t-wh_x\tau_w.
\label{psi-form}
\end{equation}

On the other hand, when \eqref{phi-u} and  \eqref{phi-w} are used in  \eqref{E2},
\begin{equation}
\phi+2w\bigl(h_x\tau_u+\tau_t+2u\tau_x-\xi_x\bigr)=0
\end{equation}
is obtained. Thus we have
\begin{equation}
\phi=2w\bigl(\xi_x-\tau_t-2u\tau_x-h_x\tau_u\bigr).
\label{phi-form}
\end{equation}

When \eqref{psi-form} and \eqref{phi-form} are used in \eqref{E4s}
\begin{equation}
\begin{aligned}
w\bigl(&h_xw\tau_{ww}-2h_x\tau_{uu}+2h_x\tau_w
+u^2\tau_{xw}+u\tau_{tw}-4u\tau_{xu}\\
&-u\xi_{xw}+w\tau_{xw}-3\tau_x-2\tau_{tu}
-\xi_{tw}+2\xi_{xu}\bigr)=0
\end{aligned}
\end{equation}
is obtained. Equations \eqref{xi-u}, \eqref{xi-w}, and their differential
consequences give
\begin{align}
\xi_{xu}&=u\tau_{xu}-w\tau_{xw},\\
\xi_{xw}&=-\tau_{xu}+u\tau_{xw},\\
\xi_{tw}&=-\tau_{tu}+u\tau_{tw}.
\end{align}
Also from \eqref{T1} it follows that
$\tau_{uu}=w\tau_{ww}+2\tau_w$. Thus
\begin{equation}
\tau_{tu}+u\tau_{xu}+w\tau_{xw}+3\tau_x+h_x\tau_{uu}=0
\label{T2}
\end{equation}
is obtained.

If the same expressions are used in \eqref{E10s}, in the first stage
\begin{equation}
h_xw\tau_{uw}-u^2\tau_{xu}+4uw\tau_{xw}-u\tau_{tu}
+u\xi_{xu}+2w\tau_{tw}
-w\tau_{xu}-2w\xi_{xw}+\xi_{tu}=0
\end{equation}
is obtained. If one uses
\begin{equation}
\xi_{tu}=u\tau_{tu}-w\tau_{tw}
\end{equation}
and  the other differential results above, it is found that
\begin{equation}
\tau_{tw}+\tau_{xu}+u\tau_{xw}+h_x\tau_{uw}=0.
\label{T3}
\end{equation}
The derivative of \eqref{T2} with respect to $w$ is
\begin{equation}
\tau_{tuw}+u\tau_{xuw}+w\tau_{xww}+4\tau_{xw}
+h_x\tau_{uuw}=0
\label{T2w}
\end{equation}
and the derivative of \eqref{T3} with respect to $u$ is
\begin{equation}
\tau_{twu}+\tau_{xuu}+\tau_{xw}+u\tau_{xuw}
+h_x\tau_{uuw}=0
\label{T3u}
\end{equation}
Difference of \eqref{T2w} and \eqref{T3u} gives
\begin{equation}
w\tau_{xww}+3\tau_{xw}-\tau_{xuu}=0.
\label{compat-a}
\end{equation}
On the other hand, from the derivative of \eqref{T1} with respect to $x$ we find
\begin{equation}
\tau_{xuu}-w\tau_{xww}-2\tau_{xw}=0.
\label{compat-b}
\end{equation}
Therefore, \eqref{compat-a} and \eqref{compat-b} give
\begin{equation}
\tau_{xw}=0,\qquad \tau_{xuu}=0
\end{equation}
Hence
\begin{equation}
\tau_x=A(t,x)u+B(t,x).
\end{equation}
We integrate this with respect to $x$ and obtain
\begin{equation}
\tau=a(t,x)u+b(t,x)+\rho(t,w,u).
\label{tau-rho}
\end{equation}
\eqref{T1}, \eqref{T2} and \eqref{T3} reduce to
\begin{align}
\rho_{uu}-w\rho_{ww}-2\rho_w&=0,\label{rho1}\\
\rho_{tw}+a_x+h_x\rho_{uw}&=0,\label{rho2}\\
a_t+\rho_{tu}+4ua_x+3b_x+h_x\rho_{uu}&=0
\label{rho3}
\end{align}
respectively.
Now, differentiating \eqref{rho1} with respect to $t$ gives
\begin{equation}
\rho_{tuu}-w\rho_{tww}-2\rho_{tw}=0.
\label{rho1t}
\end{equation}
From \eqref{rho2}
\begin{equation}\label{eq31}
\rho_{tw}=-a_x-h_x\rho_{uw}
\end{equation}
and from the $w$ derivative of \eqref{rho2} we obtain
\begin{equation}\label{eq32}
\rho_{tww}+h_x\rho_{uww}=0.
\end{equation}
If \eqref{eq31} and \eqref{eq32} are used in \eqref{rho1t}
\begin{equation}\label{eq33}
\rho_{tuu}+wh_x\rho_{uww}+2a_x+2h_x\rho_{uw}=0
\end{equation}
is obtained. Differentiating \eqref{rho1} with respect to  $u$ we have
\begin{equation}
\rho_{uuu}-w\rho_{uww}-2\rho_{uw}=0
\end{equation}
When we use this in \eqref{eq33}
\begin{equation}
\rho_{tuu}+h_x\rho_{uuu}+2a_x=0.
\label{ax-a}
\end{equation}
On the other hand, the derivative of \eqref{rho3} with respect to  $u$ gives
\begin{equation}
\rho_{tuu}+h_x\rho_{uuu}+4a_x=0.
\label{ax-b}
\end{equation}
Taking the difference of \eqref{ax-a} and \eqref{ax-b} gives
\begin{equation}
a_x=0.
\end{equation}
Since $a=a(t)$, the term $a(t)u$ can be absorbed into $\rho(t,w,u)$.
Therefore, after renaming the resulting function,
\begin{equation}
\tau=b(t,x)+R(t,w,u).
\label{tau-bR}
\end{equation}
The remaining equations are
\begin{align}
R_{uu}-wR_{ww}-2R_w&=0,\label{R1}\\
R_{tw}+h_xR_{uw}&=0,\label{R2}\\
R_{tu}+3b_x+h_xR_{uu}&=0.
\label{R3}
\end{align}
In particular, there is no $b_t$ term in the last equation.

Equations \eqref{R2} and \eqref{R3} imply
\begin{equation}
\frac{\partial}{\partial w}(R_t+h_xR_u)=0,
\qquad
\frac{\partial}{\partial u}(R_t+h_xR_u)=-3b_x.
\end{equation}
Therefore,
\begin{equation}
R_t+h_xR_u=-3b_xu+c(t,x).
\label{transport}
\end{equation}
Taking the derivative of this equation with respect to $x$ gives
\begin{equation}
h_{xx}R_u=-3b_{xx}u+c_x.
\label{xtransport}
\end{equation}
Taking derivatives with respect to $u$ and $w$, respectively, gives
\begin{align}
h_{xx}R_{uu}&=-3b_{xx},\\
h_{xx}R_{uuu}&=0,\label{branch-u}\\
h_{xx}R_{uw}&=0.\label{branch-w}
\end{align}
Thus, the branches $h_{xx}\ne0$ and $h_{xx}=0$ emerge.

\subsection{The branch $h_{xx}\ne0$}

From \eqref{branch-u}--\eqref{branch-w},
\begin{equation}
R_{uuu}=0,\qquad R_{uw}=0
\end{equation}
are found. Accordingly,
\begin{equation}
R=\frac12A(t)u^2+B(t)u+C(t,w).
\end{equation}
If this expression is used in \eqref{R1},
\begin{equation}
wC_{ww}+2C_w=A(t)
\end{equation}
is obtained. Two integrations with respect to $w$ give
\begin{equation}
R=\frac{A(t)}2(u^2+w)+B(t)u+\frac{K(t)}w+L(t)
\label{R-curved}
\end{equation}
The term $K(t)/w$ is understood locally in the region $w\ne0$. When \eqref{R-curved} is used in \eqref{transport}, the coefficients of the
powers of $u$ and $w$ give
\begin{equation}
A'=0,\qquad K'=0,\qquad B'+A h_x=-3b_x.
\end{equation}
Writing $A=A_0$ and $K=K_0$ gives
\begin{equation}
b(t,x)=-\frac{x}{3}B'(t)-\frac{A_0}{3}h(x)+d(t).
\end{equation}
Thus,
\begin{equation}
\tau=\frac{A_0}{2}(u^2+w)-\frac{A_0}{3}h+B(t)u
-\frac{x}{3}B'(t)+\frac{K_0}{w}+M(t).
\label{tau-curved}
\end{equation}
Integrating \eqref{xi-u}--\eqref{xi-w} gives
\begin{equation}
\xi=A_0\left(\frac{u^3}{3}-\frac{wu}{2}\right)
+B(t)\left(\frac{u^2}{2}-w\right)+\frac{K_0u}{w}+C(t,x).
\label{xi-curved}
\end{equation}
From \eqref{phi-form} and \eqref{psi-form},
\begin{align}
\phi&=\frac{2w}{3}\left(-A_0uh_x-uB'+xB''-3Bh_x
+3C_x-3M'\right),\label{phi-curved}\\
\psi&=\frac{A_0}{3}u^2h_x-\frac{A_0}{6}wh_x+\frac{K_0}{w}h_x
-\frac16u^2B'+\frac13xuB''+uC_x-uM'
-\frac23wB'+C_t.\label{psi-curved}
\end{align}

Now \eqref{tau-curved}--\eqref{psi-curved} are substituted into the remaining
equations \eqref{E1} and \eqref{E7}. The coefficients of $u^2$ and $w^{-1}$
in \eqref{E1} give
\begin{equation}
A_0h_{xx}=0,\qquad K_0h_{xx}=0.
\end{equation}
Since $h_{xx}\ne0$ in this branch,
\begin{equation}
A_0=K_0=0.
\end{equation}
The remaining coefficient equations are
\begin{align}
B''+9C_{xx}-6Bh_{xx}&=0,\label{cur1}\\
2xB'''-7B'h_x-6M''+9C_{tx}&=0,\label{cur2}\\
B''+6C_{xx}-3Bh_{xx}&=0,\label{cur3}\\
2C_{xx}-Bh_{xx}&=0,\label{cur4}\\
xB'''-3B'h_x-3M''+6C_{tx}&=0,\label{cur5}\\
2xB''h_x-3B(h_x)^2-3Ch_{xx}+3C_{tt}
+3C_xh_x-6M'h_x&=0.\label{cur6}
\end{align}
Equation \eqref{cur4} gives $C_{xx}=Bh_{xx}/2$. When this is used in \eqref{cur3},
$B''=0$ is obtained. Finally, when these two results are used in \eqref{cur1},
\begin{equation}
-\frac32Bh_{xx}=0
\end{equation}
is obtained. Therefore,
\begin{equation}
B=0,\qquad C_{xx}=0.
\end{equation}
\eqref{cur2} and \eqref{cur5} now become
\begin{equation}
-6M''+9C_{tx}=0,\qquad -3M''+6C_{tx}=0.
\end{equation}
Hence,
\begin{equation}
M''=0,\qquad C_{tx}=0.
\end{equation}
Thus
\begin{equation}
M=c_1t+c_0,\qquad C=c_2x+q(t).
\end{equation}
Substitution into the last equation, \eqref{cur6}, gives
\begin{equation}
[c_2 x+q(t)]h''(x)+(2c_1-c_2)h'(x)-q''(t)=0
\label{classifying}
\end{equation}
which is the classification equation. As a result, in the branch $h_{xx}\ne0$,
\begin{align}
\tau&=c_1t+c_0,\\
\xi&=c_2x+q(t),\\
\phi&=2(c_2-c_1)w,\\
\psi&=(c_2-c_1)u+q'(t),
\end{align}
and the condition \eqref{classifying} is obtained for $h$.

\subsection{The branch $h_{xx}=0$}

In this branch,
\begin{equation}
h(x)=\alpha x+\beta,\qquad h_x=\alpha
\end{equation}
where $\alpha,\beta$ are constants. Since the left-hand side of
\eqref{transport} is independent of $x$, differentiation with respect to $x$ gives
\begin{equation}
-3b_{xx}u+c_x=0.
\end{equation}
Therefore,
\begin{equation}
b_{xx}=0,\qquad c_x=0,
\end{equation}
and
\begin{equation}
b(t,x)=p(t)x+q(t),\qquad c=c(t)
\end{equation}
can be written.

Let $p(t)=P'(t)$ and define
\begin{equation}
z=u-\alpha t.
\end{equation}
Then
\begin{equation}
(\partial_t+\alpha\partial_u)z=0.
\end{equation}
The general solution of \eqref{transport}, with a function $T(t)$, can be
written as
\begin{equation}
\tau=F(w,z)+P'(t)x-3P(t)z+T(t).
\label{tau-linear}
\end{equation}
Substitution into \eqref{R1} gives
\begin{equation}
F_{zz}-wF_{ww}-2F_w=0.
\label{F-eqn}
\end{equation}
To integrate \eqref{xi-u}--\eqref{xi-w}, define the function $G(w,z)$ by
\begin{align}
G_z&=zF_z-wF_w,\label{G1}\\
G_w&=-F_z+zF_w.\label{G2}
\end{align}
The compatibility condition of these two equations is exactly
\eqref{F-eqn}. Thus,
\begin{equation}
\xi=G(w,z)+\alpha tF(w,z)-\frac32P(t)u^2+3P(t)w+D(t,x).
\label{xi-linear}
\end{equation}

Equations \eqref{tau-linear} and \eqref{xi-linear} are substituted into
\eqref{phi-form} and \eqref{psi-form}. The derivatives of $G$ are eliminated
using \eqref{G1}--\eqref{G2}. The result is
\begin{align}
\phi&=2w\left[(u-3\alpha t)P'-xP''+D_x-T'\right],
\label{phi-linear}\\
\psi&=\alpha F-3\alpha uP
+\left(\frac{u^2}{2}-3\alpha tu+2w\right)P'
-uxP''+uD_x-uT'+D_t.
\label{psi-linear}
\end{align}

Finally, \eqref{tau-linear}, \eqref{xi-linear}, \eqref{phi-linear}, and
\eqref{psi-linear} are substituted into \eqref{E1} and \eqref{E7}.
Separating the coefficients with respect to $z$ and $w$ gives
\begin{align}
D_{xx}&=0,\label{lin1}\\
P''&=0,\label{lin2}\\
5\alpha P'+2T''-3D_{tx}&=0,\label{lin3}\\
-3\alpha P'-T''+2D_{tx}&=0.\label{lin4}
\end{align}
From \eqref{lin3}--\eqref{lin4},
\begin{equation}
T''=-\alpha P',\qquad D_{tx}=\alpha P'.
\end{equation}
Since $P''=0$,
\begin{equation}
P(t)=p_1t+p_0
\end{equation}
and
\begin{equation}
T(t)=-\frac{\alpha p_1}{2}t^2+t_1t+t_0.
\end{equation}
Equations $D_{xx}=0$ and $D_{tx}=\alpha p_1$ give
\begin{equation}
D(t,x)=\left[\alpha(p_1t+p_0)+c_1\right]x+d_0(t).
\end{equation}
The remaining constant term in \eqref{E7} gives
\begin{equation}
d_0''=6\alpha^2p_1t+2\alpha^2p_0-\alpha c_1+2\alpha t_1,
\end{equation}
hence
\begin{equation}
\begin{aligned}
d_0(t)={}&\alpha^2p_1t^3
+\left(\alpha^2p_0-\frac{\alpha c_1}{2}+\alpha t_1\right)t^2
+c_2t+c_3.
\end{aligned}
\end{equation}

As a result, the general infinitesimal coefficients in the $h_{xx}=0$ branch are
\begin{align}
\tau={}&F(w,z)+p_1x-3(p_1t+p_0)z
-\frac{\alpha p_1}{2}t^2+t_1t+t_0,
\label{final-tau}\\
\xi={}&G(w,z)+\alpha tF(w,z)-\frac32(p_1t+p_0)u^2
+3(p_1t+p_0)w\notag\\
&+\left[\alpha(p_1t+p_0)+c_1\right]x+\alpha^2p_1t^3\notag\\
&+\left(\alpha^2p_0-\frac{\alpha c_1}{2}+\alpha t_1\right)t^2
+c_2t+c_3,
\label{final-xi}\\
\phi={}&2w\left[p_1z+\alpha p_0+c_1-t_1\right],
\label{final-phi}\\
\psi={}&\alpha F(w,z)
+p_1\left(\frac{z^2}{2}+2w+\alpha x-3\alpha tz
-\frac{\alpha^2t^2}{2}\right)\notag\\
&+(c_1-2\alpha p_0)z+t_1(\alpha t-z)+c_2.
\label{final-psi}
\end{align}
Here,
\begin{equation}
z=u-\alpha t,
\end{equation}
$p_1,p_0,t_1,t_0,c_1,c_2,c_3$ are arbitrary constants, and $F$ and $G$
satisfy \eqref{F-eqn}, \eqref{G1}, and \eqref{G2}.

When all finite-dimensional parameters are set to zero, the
infinite-dimensional hodograph family
\begin{equation}
X_{F,G}=F(w,z)\partial_t+\bigl(G(w,z)+\alpha tF(w,z)\bigr)\partial_x
+\alpha F(w,z)\partial_u
\end{equation}
is obtained.

\end{document}